\documentclass[leqno,12pt]{amsart}

\usepackage{amssymb}
\usepackage{amscd}              
\usepackage[all]{xy}            
\xyoption{curve}
\CompileMatrices

\newcommand{\Aa}{{\mathbb{A}}}

\newcommand{\bB}{{\bf B}}

\newcommand{\bG}{{\bf G}}

\newcommand{\bU}{{\bf U}}

\newcommand{\comp}{\circ}

\newcommand{\diag}{\text{diag}}        
\newcommand{\divi}{\text{div}\, }                
\newcommand{\dual}{\vee}
\newcommand{\eins}{{\mathbb{I}}}

\newcommand{\E}{{\mathcal{E}}}

\newcommand{\F}{{\mathcal{F}}}
\newcommand{\Fl}{\text{\bf Fl}}

\newcommand{\GL}{\text{GL}}
\newcommand{\GLn}{\text{GL}_n}

\newcommand{\Hom}{\text{Hom}}

\newcommand{\isomorph}{\cong}           
\newcommand{\isomto}{\overset{\sim}{\rightarrow}}

\newcommand{\Isom}{\text{Isom}}

\newcommand{\KGL}{\text{KGL}}
\newcommand{\KGLn}{\text{KGL}_n}
\newcommand{\KT}{{\text{KT}}}
\newcommand{\M}{{\mathcal{M}}}

\newcommand{\Ob}{\overline{\text{\bf O}}}
\newcommand{\Oo}{{\mathcal{O}}}                  
\newcommand{\Oplus}{\bigoplus}

\newcommand{\PGL}{\text{PGL}}
\newcommand{\PGLb}{\overline{\PGL}}
\newcommand{\Pic}{\text{Pic}}

\newcommand{\Pp}{{\mathbb{P}}}                   

\newcommand{\Proj}{\text{Proj}\,}
\newcommand{\Q}{{\mathbb{Q}}}                    

\newcommand{\rk}{\text{rk}\, }   

\newcommand{\s}{{\sf s}}

\newcommand{\tensor}{\otimes}

\newcommand{\Tensor}{\bigotimes}

\newcommand{\To}{\longrightarrow}
\newcommand{\Tt}{\tilde{T}}

\newcommand{\univ}{\text{univ}}

\newcommand{\Xb}{\overline{X}}

\newcommand{\Z}{{\mathbb{Z}}}                    

\newtheorem{theorem}{Theorem}[section]
\newtheorem{proposition}[theorem]{Proposition}
\newtheorem{lemma}[theorem]{Lemma}
\newtheorem{corollary}[theorem]{Corollary}

\theoremstyle{definition}
\newtheorem{definition}[theorem]{Definition}
\newtheorem{remark}[theorem]{Remark}

\makeatletter
\def\overunderbraces #1#2#3{{%
 \baselineskip\z@skip \lineskip4\p@ \lineskiplimit4\p@
 \displaystyle  
 \setbox\z@\vbox{\ialign{&\hfil${}##{}$\hfil\cr
   \global\let\br\br@label #1\cr 
   \global\let\br\br@down #1\cr   
   #2\cr 
 }}
 \dimen@-\ht\z@ 
 \setbox\z@\vbox{\ialign{&\hfil${}##{}$\hfil\cr
   \global\let\br\br@label #1\cr 
   \global\let\br\br@down #1\cr   
   #2\cr 
   \global\let\br\br@up #3\cr 
   \global\let\br\br@label #3\cr   
 }}
 \advance\dimen@\ht\z@ 
 \lower\dimen@\hbox{\box\z@} 
}}
\def\br@up#1#2{\multispan{#1}\upbracefill}
\def\br@down#1#2{\multispan{#1}\downbracefill}
\def\br@label#1#2{\multispan{#1}\hidewidth $#2$\hidewidth}
\makeatother

\begin{document}

\title[Global Sections of Line Bundles on $\KGL(E,F)$]
{Global Sections of Line Bundles on a Wonderful Compactification
 of the General Linear Group}
\author[Ivan Kausz]{Ivan Kausz}
\thanks{Partially supported by the DFG} 
\date{September 5, 2004}
\address{NWF I - Mathematik, Universit\"{a}t Regensburg, 93040 Regensburg, 
Germany}
\email{ivan.kausz@mathematik.uni-regensburg.de}
\begin{abstract}
In a previous paper \cite{kgl} I have constructed a compactification
       $\KGL_n$ of the general linear group $\GL_n$, 
       which in many respects is analogous
       to the so called wonderful compactification of adjoint semisimple
       algebraic groups as studied by De Concini and Procesi.
       In particular there is an action of $\bG=\GL_n\times\GL_n$ on this
       compactification.
       In this paper we show how the space of global section of an
       arbitrary $\bG$-linearized line bundle on $\KGL_n$ and its
       orbit-closures decomposes
       into a direct sum of simple $\bG$-modules.
\end{abstract}

\maketitle

\tableofcontents


\section{Introduction}

Let $k$ be a field of characteristic zero and let $E$ and $F$ be two 
$n$-dimensional vector spaces.
In \cite{kgl} we have introduced a certain compactification $\KGL(E,F)$ 
of the variety $\Isom(E,F)\isomorph \GL_n$ of linear isomorphisms from 
$E$ to $F$ which in many respects
is analogous to De Concinis and Procesis so called 
{\em wonderful compactification} 
of adjoint semi-simple algebraic groups (cf. \cite{CP}).

In particular, there is a natural action of the group 
$\bG:=\GL(E)\times\GL(F)$ on $\KGL(E,F)$ extending the one arising
from right and left multiplication on $\Isom(E,F)$.
Furthermore $\KGL(E,F)$ is smooth, the boundary,
i.e. the complement  of $\Isom(E,F)$ in  $\KGL(E,F)$,
is a divisor with normal crossings 
and the closures of the orbits of the $\bG$-action are precisely the 
nonempty intersections of the irreducible components of the boundary. 

We will see in \ref{picard} below that the Picard group of 
$\KGL(E,F)$ is generated
by (the ideal sheaves of) the boundary components
$Z_0,\dots,Z_{n-1}$ and $Y_0,\dots,Y_{n-1}$. Every line bundle 
expressed in terms of these generators is equipped with a canonical
linearization of the $\bG$-action and thus the space of global sections
of its restriction to some orbit closure is naturally a finite dimensional 
$\bG$-module.

In this paper we show how such a space of global sections decomposes into
a direct sum of simple $\bG$-modules. More precisely, we prove the
following 

\vspace{2mm}
\noindent
{\bf Theorem:}
{\em
(Cf. Theorem \ref{theorem} for the exact formulation). 
Let $L$ be a $\bG$-linearized line bundle 
of the form $L=\Oo(\sum(m_iZ_i+l_iY_i))$ on $\KGL(E,F)$ and
let $I,J\subseteq[0,n-1]$ be subsets such that the intersection
$
\Ob_{IJ}=(\cap_{i\in I}Z_i)\cap(\cap_{j\in J}Y_j)
$ 
is nonempty.
Then the decomposition of the $\bG$-module $H^0(\Ob_{IJ},i_{IJ}^*L)$
into simple submodules is given by a canonical isomorphism
$$
H^0(\Ob_{IJ},i_{IJ}^*L)\isomto 
\Oplus_{(a,b)\in A_{IJ}(L)}H^0(\Fl,\Oo_{\Fl}(a,b))
\quad,
$$
where $A_{IJ}(L)\subset\Z^n\times\Z^n$ is a finite set defined 
explicitly in terms of
$I$, $J$ and $L$, where $\Fl$ is the product of the two complete 
flag manifolds associated to the vector spaces $E$ and $F$ respectively
and where $\Oo_{\Fl}(a,b)$ is the product specified by $(a,b)$ 
of successive quotients of tautological vector bundles on $\Fl$.
}
\vspace{2mm}

In \cite{degeneration} we have shown the relevance of $\KGL(E,F)$ for 
the Gieseker type degeneration of moduli stacks of vector bundles on
curves: The normalization of the moduli stack of Gieseker vector bundles
on an irreducible nodal curve with one singularity is isomorphic to
$\KGL(\E,\F)$, where $\E$ and $\F$ are certain vector bundles on the moduli
stack of vector bundles on the normalization of the curve.
In a forthcoming paper we will apply the
results of the present paper to obtain a canonical decomposition of
generalized theta functions on the moduli stack of Gieseker vector bundles
(cf. \cite{factorisation}).

Our proof of Theorem \ref{theorem} is 
inspired by \cite{CP} \S 8, where the cohomology of line bundles 
on complete symmetric varieties is computed. 
At one notable point however we have to argue differently, since
to show that certain simple submodules occur in the space of
global sections De Concini and Procesi make use of the fact that
certain line bundles are ample (cf. \cite{CP}, Proposition 8.4), 
and it  turns out (cf. \ref{not ample}) that the corresponding 
statement is false in the case of $\KGL(E,F)$. Instead, we produce 
(in \ref{divisor}) explicit sections which 
generate the simple submodules in question.

After finishing this paper I have learned that A. Tchoudjem
has studied the cohomology of line bundles on compactifications
of arbitrary reductive groups \cite{Tchoudjem}. Part of our
result can probably be deduced from his, but certainly not all, 
since  he does not deal with cohomology of the strata and does
not obtain a {\em canonical} decomposition.

This paper has been written during a stay at the Tata Institute of
Fundamental Research in Bombay. Its hospitality is gratefully 
acknowledged.

\section{Notation}
If $p\leq q$ are two integers, we denote
by $[p,q]$ the set of all integers $i$ with $p\leq i\leq q$.

\section{Preliminary results}
\label{preliminary}

In this chapter we recall some results from \cite{kgl}
which we will need in the following chapters.

Let $k$ be a field of characteristic zero.
We fix two $k$-vector spaces $E$ and $F$ of rank $n$.
In \cite{kgl} I have defined a compactification $\KGL(E,F)$ 
of the scheme $\Isom(E,F)$ by the following construction:
Let $X^{(0)}:=\Pp(\Hom(E,F)^\dual\oplus k)$ and define
for $i=0, \dots n-1$ the closed subschemes 
\begin{eqnarray*}
Y_i^{(0)}
&:=&
\{(f:a)\in X^{(0)}\ |\
\rk(f)\leq i\}
\\
Z_{i}^{(0)}
&:=&
\{(f:0)\in X^{(0)}\ |\
\rk(f)\leq n-i\}
\end{eqnarray*}
of $X^{(0)}$. 
In other words, after choosing a basis for $E$ and for
$F$ we can identify the scheme 
$
X^{(0)}
$
with
$
\Proj(k[x_{00},x_{ij}(i,j\in[1,n])])
$ and
the subscheme $Y_p^{(0)}$ (the subscheme $Z_{n-p}^{(0)}$)
belongs to the
homogenous ideal generated by the $(p+1)\times(p+1)$-subminors
of the matrix $(x_{ij})_{i,j}$ (by these minors and $x_{00}$). 
These subschemes satisfy inclusion relations as indicated
in the following diagram:
$$
\xymatrix@R-1.2pc{
Y_0^{(0)} \ar@<-.5ex>@{^{(}->}[r] & 
Y_1^{(0)} \ar@<-.5ex>@{^{(}->}[r] & 
\quad\dots\quad \ar@<-.5ex>@{^{(}->}[r]     & 
Y_{n-1}^{(0)}             & 
\\
                              & 
Z_{n-1}^{(0)} \ar@<-.5ex>@{^{(}->}[r] \ar@{^{(}->}[u]& 
\quad\dots\quad \ar@<-.5ex>@{^{(}->}[r]         & 
Z_1^{(0)} \ar@<-.5ex>@{^{(}->}[r] \ar@{^{(}->}[u]    & 
Z_0^{(0)} }
$$
By definition, $\KGL(E,F)$ is the result of successively
blowing up the scheme $X^{(0)}$ as follows:
$$
\xymatrix{
X^{(0)} &
X^{(1)} \ar[l] &
X^{(2)} \ar[l] &
\dots\dots \ar[l] &
X^{(n-1)}=\KGL(E,F) \ar[l] 
} 
$$
Here, $X^{(1)}$ is the result of blowing up
$X^{(0)}$ in the (disjoint) union of the subschemes $Y^{(0)}_0$
and $Z^{(0)}_{n-1}$. 
Generally, in the $i$-th step we define
$Y^{(i)}_{j-1},Z^{(i)}_{n-j}\subset X^{(i)}$
to be the  proper transforms of $Y^{(i-1)}_{j-1}$ and $Z^{(i-1)}_{n-j}$
respectively if $j\neq i$ and to be the exceptional divisors lying
above $Y_{i-1}^{(i-1)}$ and $Z_{n-i}^{(i-1)}$ respectively if $j=i$.
Then it turns out that the subschemes $Y^{(i)}_i$ and $Z^{(i)}_{n-i-1}$
are smooth and disjoint and thus the blowing up of $X^{(i)}$ in 
$Y^{(i)}_i\cup Z^{(i)}_{n-i-1}$ is a smooth projective 
variety $X^{(i+1)}$.

We have a natural open embedding $\Isom(E,F)\subset X^{(0)}$.
Since the centers of blowing up are in the complement of $\Isom(E,F)$,
we can regard $\Isom(E,F)$ as an open subset of $\KGL(E,F)$.
By \cite{kgl}, 4.2 the complement of $\Isom(E,F)$
in $\KGL(E,F)$ is a divisor with normal crossings
whose irreducible components are $Z_1,\dots,Z_{n-1}$ and $Y_1,\dots,Y_{n-1}$,
where $Z_i:=Z_i^{(n-1)}$ and $Y_i:=Y_i^{(n-1)}$.

After the choice of a basis for $E$ and $F$ there are canonical
rational functions 
\begin{eqnarray*}
y_{ji}, z_{ij} & & (1\leq i<j\leq n) \\
t_i/t_0 & & (1\leq i\leq n)
\end{eqnarray*}
on $\KGL(E,F)$
which are related to the coordinate functions $x_{ij}/x_{00}$
on $X^{(0)}$
by the matrix equation
$$
\label{matrix decomposition}
\left[\frac{x_{ij}}{x_{00}}\right] = 
\left[
\vcenter{
\hbox{
\xymatrix @R-2.5pc@C-2.5pc@!{
1 \ar@{.}[rrrrrddddd]
      &     &     &     &      & 
\ar@{--} '[lllll]+<1.3ex,0ex> '[ddddd]+<0ex,1.3ex> []+0 \\
      &     &     &     &  0   &     \\
      &     &     &     &      &     \\
      & y_{ij}&     &     &      &     \\
      &     &     &     &      &     \\
\ar@{--} '[uuuuu]+<0ex,-1.3ex> '[rrrrr]+<-1.3ex,0ex> []+0
      &     &     &     &      &  1 }}}
\right]
\left[
\vcenter{
\hbox{
\xymatrix @R-3.7pc@C-3.7pc@! {
t_1/t_0 \ar@{.}[rrrrrddddd]
      &     &     &     &     &      \\
      &     &     &  0  &     &      \\
      &     &     &     &     &      \\
      &     &     &     &     &      \\
      &     &  0  &     &     &      \\
      &     &     &     &     &  t_n/t_0 }}}
\right]
\left[
\vcenter{
\hbox{
\xymatrix @R-2.5pc@C-2.5pc@!{
1 \ar@{.}[rrrrrddddd]
      &     &     &     &      & 
\ar@{--} '[lllll]+<1.3ex,0ex> '[ddddd]+<0ex,1.3ex> []+0 \\
      &     &     &     &  z_{ij} &     \\
      &     &     &     &      &     \\
      & 0   &     &     &      &     \\
      &     &     &     &      &     \\
\ar@{--} '[uuuuu]+<0ex,-1.3ex> '[rrrrr]+<-1.3ex,0ex> []+0
      &     &     &     &      &  1 }}}
\right]
$$

\begin{proposition}
\label{X(ell)}
Fix a basis for $E$ and for $F$.
Let $\ell\in[0,n]$.
Let $\iota_{\ell}:[1,n+1]\to[0,n]$ be the bijection
such that
$\iota_{\ell}(\ell+1)=0$ and such that it induces an
increasing map from the complement of $\ell+1$ onto $[1,n]$.

There is an open subscheme 
$X(\ell)\subset\KGL(E,F)$ which is isomorphic to
the $n^2$-dimensional affine space, such that the coordinate
functions on $X(\ell)$ are the restrictions of the rational
functions $y_{ji}, z_{ij}$ ($1\leq i<j\leq n$) and the
rational functions
$$
t_{\iota_{\ell}(i+1)}/t_{\iota_{\ell}(i)}
\quad
(1\leq i\leq n)
\quad.
$$

Furthermore, the intersection of $Y_i$ (of $Z_{n-i-1}$) 
with $X(\ell)$ is empty if $0\leq i\leq\ell-1$
(if $\ell\leq i\leq n-1$),
and the equation on $X(\ell)$ of the divisor
$Y_i$ (of the divisor $Z_{n-i-1}$)
is $t_{\iota_{\ell}(i+2)}/t_{\iota_{\ell}(i+1)}$
if $\ell\leq i\leq n-1$ (if $0\leq i\leq\ell-1$).
\end{proposition}

\begin{proof}
This is a special case of \cite{kgl}, Proposition 4.1.
\end{proof}

In \cite{kgl}, 5.5 I have shown that $\KGL(E,F)$ can be regarded as 
a moduli space for certain diagrams of vector bundles.
To formulate the precise statement, we have to introduce
some definitions.

Let $T$ be a $k$-scheme and let $\E$, $\F$ be two locally free
$\Oo_T$-modules of rank $n$. A 
{\em bf-morphism} from $\E$ to $\F$ is a tuple
$
\gamma=(L, \lambda, \E\to\F, \F\to L\tensor\E, r)
$
where $L$ is an invertible $\Oo_T$-module, $\lambda$ is a section of $L$,
the arrows $\E\to\F$ and $\F\to L\tensor\E$ are $\Oo_T$-module morphisms and 
$r$ is an integer between $0$ and $n$ 
such that locally on $T$ there exist isomorphisms
\begin{eqnarray*}
\E &\isomto& r\Oo_T\oplus(n-r)\Oo_T \\
\F &\isomto& r\Oo_T\oplus(n-r)L
\end{eqnarray*}
with the property that via these isomorphisms the morphisms
$
\E\to\F
$
and 
$\F\to L\tensor\E$
are expressed by the diagonal matrices
$$
\text{
$
\left[
\begin{array}{cc}
\eins_r & 0 \\
0 & \lambda\eins_{n-r}
\end{array}
\right]
$
\quad and \quad
$
\left[
\begin{array}{cc}
\lambda\eins_r & 0 \\
0 & \eins_{n-r}
\end{array}
\right]
$
}
$$
respectively.
We will often use the following more suggestive notation for the 
bf-morphism $\gamma$:
$$
\gamma=\left(
\xymatrix@C2.5ex{
\E
\ar[rr]^r_{(L,\lambda)}
& & 
\F
\ar@/_1.2pc/|{\tensor}[ll]
}
\right)
\quad.
$$

Let $T$, $\E$, $\F$ be as above.
A {\em generalized isomorphism} $\Phi$ from $\E$ to $\F$ is a sequence of
bf-morphisms connected as follows:
$$
\xymatrix@C=1.2ex{
\E
\ar@/^1.2pc/|{\tensor}[rr]
& &
E_1 
\ar[ll]_0^{(M_0,\mu_0)}
\ar@/^1.2pc/|{\tensor}[rr]
& &
E_2
\ar[ll]_1^{(M_1,\mu_1)}
& 
\dots
& 
E_{n-1}
\ar@/^1.2pc/|{\tensor}[rrr]
& & &
E_n
\ar[lll]_{n-1}^{(M_{n-1},\mu_{n-1})}
\ar[r]^\sim
& 
F_n
\ar[rrr]^{n-1}_{(L_{n-1},\lambda_{n-1})}
& & &
F_{n-1}
\ar@/_1.2pc/|{\tensor}[lll]
& 
\dots
&
F_2
\ar[rr]^1_{(L_1,\lambda_1)}
& &
F_1
\ar[rr]^0_{(L_0,\lambda_0)}
\ar@/_1.2pc/|{\tensor}[ll]
& &
\F
\ar@/_1.2pc/|{\tensor}[ll]
}
$$
which has properties for which we refer the reader to
\cite{kgl} 5.2, since they will not be of importance here.
Two generalized isomorphisms $\Phi$ (as above) and $\Phi'$ 
(with primed ingredients)
from $\E$ to $\F$ are said to be {\em equivalent}, if there
exist isomorphisms 
$E_i\isomto E'_i$, 
$F_i\isomto F'_i$, 
$M_i\isomto M'_i$, 
$L_i\isomto L'_i$,
such that all the obvious diagrams commute.
Theorem 5.5 in \cite{kgl} can now be formulated as 
follows:

\begin{theorem}
\label{Phi_univ}
The variety $\KGL(E,F)$ represents the functor which
to a $k$-scheme $T$ associates the set of all equivalence
classes of generalized isomorphisms from 
$E\tensor\Oo_T$ to $F\tensor\Oo_T$.
In particular, there is a universal generalized isomorphism
$\Phi_{\univ}:$
$$
\xymatrix@C=1.4ex{
E\tensor\Oo
\ar@/^1.2pc/|{\tensor}[rr]
& &
E_1 
\ar[ll]_0^{(M_0,\mu_0)}
& 
\dots
& 
E_{n-1}
\ar@/^1.2pc/|{\tensor}[rrr]
& & &
E_n
\ar[lll]_{n-1}^{(M_{n-1},\mu_{n-1})}
\ar[r]^\sim
& 
F_n
\ar[rrr]^{n-1}_{(L_{n-1},\lambda_{n-1})}
& & &
F_{n-1}
\ar@/_1.2pc/|{\tensor}[lll]
& 
\dots
&
F_1
\ar[rr]^0_{(L_0,\lambda_0)}
& &
F\tensor\Oo
\ar@/_1.2pc/|{\tensor}[ll]
}
$$
The global sections $\mu_i$ and $\lambda_i$ of the line
bundles $M_i$ and $L_i$ vanish exactly along the boundary
components $Z_i$ and $Y_i$ respectively.
Thus we have canonincal isomorphisms
$M_i=\Oo(Z_i)$ and $L_i=\Oo(Y_i)$ which identify
$\mu_i$ and $\lambda_i$ with the canonical $1$-sections of $\Oo(Z_i)$
and of $\Oo(Y_i)$ respectively.  
\end{theorem}

In \cite{kgl}, \S 6 I have shown that a bf-morphism
$
\gamma=(L, \lambda, \E\to\F, \F\to L\tensor\E, p)
$
induces canonical morphisms
\begin{eqnarray*}
\wedge^r\E & \to & (L^{\dual})^{\tensor\max(0,r-p)}\tensor\wedge^r\F
\\
\wedge^r\F & \to & L^{\tensor\min(r,n-p)}\tensor\wedge^r\E
\quad,
\end{eqnarray*}
which I call the exterior powers of $\gamma$.
Given a generalized isomorphism $\Phi$:
$$
\xymatrix@C=1.2ex{
\E
\ar@/^1.2pc/|{\tensor}[rr]
& &
E_1 
\ar[ll]_0^{(M_0,\mu_0)}
\ar@/^1.2pc/|{\tensor}[rr]
& &
E_2
\ar[ll]_1^{(M_1,\mu_1)}
& 
\dots
& 
E_{n-1}
\ar@/^1.2pc/|{\tensor}[rrr]
& & &
E_n
\ar[lll]_{n-1}^{(M_{n-1},\mu_{n-1})}
\ar[r]^\sim
& 
F_n
\ar[rrr]^{n-1}_{(L_{n-1},\lambda_{n-1})}
& & &
F_{n-1}
\ar@/_1.2pc/|{\tensor}[lll]
& 
\dots
&
F_2
\ar[rr]^1_{(L_1,\lambda_1)}
& &
F_1
\ar[rr]^0_{(L_0,\lambda_0)}
\ar@/_1.2pc/|{\tensor}[ll]
& &
\F
\ar@/_1.2pc/|{\tensor}[ll]
}
$$
over a scheme $T$ we can compose the exterior powers of the
bf-morphisms occuring in it and can thus define the exterior power
$$
\wedge^r\Phi:\ \
\bigwedge^r\E
\ \ \To\ \ 
\Tensor_{\nu=1}^{r}
\left(
\Tensor_{i=0}^{\nu-1}L_i^{\dual}
\tensor
\Tensor_{i=0}^{n-\nu}M_i
\right)
\tensor
\bigwedge^r\F
$$
of $\Phi$.
If $\E=\F=n\Oo_T$ then $\wedge^r\E$ and $\wedge^r\F$ have
a natural direct sum decomposition into copies of $\Oo_T$
indexed by the subsets of cardinality $r$ of $[1,n]$.
Thus we have canonical inclusion and projection morphisms
$
\iota_A:\Oo_T \to \wedge^r\E
$,
$
\pi_B:\wedge^r\F\to\Oo_T
$
for such subsets $A,B\subseteq[1,n]$
and we can define the section
$\det_{A,B}\Phi:=\pi_A\comp(\wedge^r\Phi)\comp\iota_B$
of the line bundle 
$
\tensor_{\nu=1}^r
(
\tensor_{i=0}^{\nu-1}L_i^{\dual}
\tensor
\tensor_{i=0}^{n-\nu}M_i
)
$.
\begin{proposition}
\label{det Phi}
\begin{enumerate}
\item
Fix a basis of $E$ and $F$, let $\ell\in[0,n]$ and let $X(\ell)$ 
be the corresponding open subset of $\KGL(E,F)$ defined in
Proposition \ref{X(ell)}. 
Then $X(\ell)$ is the largest subset in the complement
of the divisors $Y_0,\dots,Y_{\ell-1}$ and
$Z_0,\dots,Z_{n-\ell-1}$ for which the sections
$\det_{[1,r],[1,r]}\Phi_{\univ}$ are nowhere vanishing for
$r=1,\dots,n$.
\item
The morphism $\wedge^n\Phi_{\univ}$
is nowhere vanishing.
\end{enumerate}
\end{proposition}

\begin{proof}
The first statement follows from \cite{kgl} 7.4 and 4.3.
The second statement is a consequence of loc. cit. 6.5.
\end{proof}

Let $T$ be a scheme, let $\E$, $\F$ be two vector
bundles of rank $n$ on $T$, and let $\gamma$
be a bf-morphism from $\E$ to $\F$.
An automorphism $g$ of $\E$ (an automorphism $h$ of $\F$)
can be composed in an obvious way with $\gamma$
to give a new bf-morphism $\gamma g$ (a new bf-morphism $h\gamma$)
from $\E$ to $\F$.
Thus, if $\Phi$ is a generalized isomorphism
from $\E$ to $\F$, we get a new generalized automorphism
$h\Phi g$ from $\E$ to $\F$ by composing the two
outer bf-morphisms in $\Phi$ with $h$ and $g$.

\begin{definition}
\label{operation}
Let $\bG$ be the product of the two algebraic
$k$-groups $\GL(E)$ and $\GL(F)$.
There is a natural operation of $\bG$ on $\KGL(E,F)$
extending the one on $\Isom(E,F)$.
In terms of $R$-valued points ($R$ a $k$-algebra) it is
given by 
$$
g\cdot\Phi=g_2\Phi g_1^{-1}
\quad,
$$
where $g=(g_1,g_2)\in\bG(R)$ and $\Phi$ is a generalized
isomorphism from $E\tensor R$ to $F\tensor R$.
\end{definition}

\begin{corollary}
\label{B-invariant}
Let $(v_1,\dots,v_n)$ and $(w_1,\dots,w_n)$ be a basis
for $E$ and $F$ respectively.
Let $B_1\subset\GL(E)$ and $B_2\subset\GL(F)$
be the Borel subgroups consisting of linear automorphisms
fixing the flags
$$
\text{
$
\{0\}\subset\langle v_1\rangle\subset
            \langle v_1,v_{2}\rangle\subset\dots\subset E
$
\qquad and \qquad
$
\{0\}\subset\langle w_n\rangle\subset
            \langle w_n,w_{n-1}\rangle\subset\dots\subset F
$
}
$$
respectively. 
For each $\ell\in[0,n]$ the open subset 
$X(\ell)$ is  invariant by the operation of 
$\bB:=B_1\times B_2\subset\bG$
on $\KGL(E,F)$ 
\end{corollary}

\begin{proof}
Let $R$ be a $k$-algebra. By \ref{Phi_univ} and 
\ref{det Phi} (1) an $R$-valued
point of $X(\ell)$ is given by generalized isomorphism
$\Phi$ from $E\tensor R$ to $F\tensor R$ such that the sections
$\det_{[1,r],[1,r]}\Phi$ ($r\in[1,n]$),
$\lambda_0,\dots,\lambda_{\ell}$
and
$\mu_0,\dots,\mu_{n-\ell}$
are nowhere vanishing.
We have to show that for each $(g_1,g_2)\in\bB(R)$ the
generalized isomorphism $\Phi'=g_2\Phi g_1^{-1}$ has again
this property.
Since the sections $\lambda_i$ and $\mu_i$ are the same
in $\Phi'$ and in $\Phi$, it suffices to show that
the quotients 
$\det_{[1,r],[1,r]}\Phi'/\det_{[1,r],[1,r]}\Phi$ ($r\in[1,n]$)
are in $R^{\times}$.
But since $g_1\in B_1(R)$
and $g_2\in B_2(R)$,
there exist $u_1,u_2\in R^{\times}$, such that
\begin{eqnarray*}
\wedge^rg_1^{-1}(v_1\wedge\dots\wedge v_r)
&=&
u_1\cdot v_1\wedge\dots\wedge v_r
\\
\pi_{[1,r]}\left(\wedge^rg_2(w_{j_1}\wedge\dots\wedge w_{j_r})\right)
&=&
\left\{
\begin{array}{ll}
u_2 & \text{, for $\{j_1,\dots,j_r\}=[1,r]$,}
\\
0 & \text{, else.}
\end{array}
\right.
\end{eqnarray*}
Therefore we have
$\det_{[1,r],[1,r]}\Phi'=u_1u_2\det_{[1,r],[1,r]}\Phi$.
\end{proof}

\begin{definition}
From now on the symbols $L_i, M_i$ will denote the line bundles
which occur in the universal generalized isomorphism
$\Phi_{\univ}:$
$$
\xymatrix@C=1.4ex{
E\tensor\Oo
\ar@/^1.2pc/|{\tensor}[rr]
& &
E_1 
\ar[ll]_0^{(M_0,\mu_0)}
& 
\dots
& 
E_{n-1}
\ar@/^1.2pc/|{\tensor}[rrr]
& & &
E_n
\ar[lll]_{n-1}^{(M_{n-1},\mu_{n-1})}
\ar[r]^\sim
& 
F_n
\ar[rrr]^{n-1}_{(L_{n-1},\lambda_{n-1})}
& & &
F_{n-1}
\ar@/_1.2pc/|{\tensor}[lll]
& 
\dots
&
F_1
\ar[rr]^0_{(L_0,\lambda_0)}
& &
F\tensor\Oo
\ar@/_1.2pc/|{\tensor}[ll]
}
$$
on $\KGL(E,F)$.
\end{definition}

From definition \ref{operation} it 
is clear that the line bundles $M_i$ and
$L_i$ are canonically $\bG$-linearized. 
Notice that also the trivial line bundles $\det(E)\tensor_k\Oo$ and 
$\det(F)\tensor_k\Oo$ carry canonical nontrivial $\bG$-linearization.

\begin{lemma}
\label{picard}
There is a canonical isomorphism of $\bG$-linearized
line bundles on $\KGL(E,F)$:
$$
(\det E)^{-1}\tensor_k\Tensor_{i=0}^{n-1}M_i^{n-i}=
(\det F)^{-1}\tensor_k\Tensor_{i=0}^{n-1}L_i^{n-i}
\quad.
\eqno(*)
$$
The Picard group of the variety $\KGL(E,F)$ 
is generated by the isomorphism classes of the line bundles
$M_i$ and $L_i$ ($i\in[0,n-1]$) and the only relations come from the
isomorphism $(*)$.
\end{lemma}

\begin{proof}
The first statement follows from \ref{det Phi}(2),
since it says that the canonical morphism
$$
\wedge^n\Phi_{\univ}:
\det(E)\tensor_k\Oo\to\
\Tensor_{i=1}^n
\left(
\Tensor_{j=0}^{i-1}L_j^{-1} \tensor
\Tensor_{j=0}^{n-i}M_j
\right)
\tensor_k\det(F)
$$
is nowhere vanishing.

Recall that $\KGL(E,F)$ is defined as the result of a successive
blowing up $X^{(i)}\to X^{(i-1)}$
along disjoint and smooth irreducible subschemes
$Y_{i-1}^{(i-1)}$ and $Z_{n-i}^{(i-1)}$ of $X^{(i-1)}$ of 
codimension $\geq 2$.
Therefore the divisor class group of $X^{(i)}$ is the direct sum
of the divisor class group of $X^{(i-1)}$ and the free abelian
group generated by the two divisors $Y_{i-1}^{(i)}$ and $Z_{n-i}^{(i)}$.
Now the divisor class group of $X^{(0)}$ is generated by the
hyper-plane $Z_0^{(0)}$; therefore by induction it follows that
the classes of $Y_0,\dots,Y_{n-2}$, $Z_0,\dots,Z_{n-1}$ 
freely generate the divisor class group of $\KGL(E,F)$.
\end{proof}

For each pair of subsets $I,J\subseteq[0,n-1]$ with $\min(I)+\min(J)\geq n$
we define the closed subscheme
$\Ob_{IJ}=\Ob_{IJ}(E,F)$ in $\KGL(E,F)$ as the intersection of  the 
components $Z_i$ ($i\in I$) and $Y_j$ ($j\in J$).
As shown in \cite{kgl} \S 9, the subschemes $\Ob_{IJ}$
are precisely the closures of the orbits of $\bG$ acting on $\KGL(E,F)$.

If $I\supseteq I'$ and $J \supseteq J'$ then we have 
$\Ob_{IJ}\subseteq\Ob_{I'J'}$.
In particular, we have $\Ob_{\emptyset\emptyset}=\KGL(E,F)$ and
the smallest of the closed subschemes $\Ob_{IJ}$ are of
the form 
$$
\Ob_{r,s}:=\Ob_{[s,n-1],[r,n-1]}
$$
for $r,s\in[0,n]$, $r+s=n$. (The set $[s,n-1]$ contains
$r$ elements while the set $[r,n-1]$ contains $s$ elements,
that's why we write $\Ob_{r,s}$ instead of $\Ob_{s,r}$).
Let
$$
\xymatrix@R=1.1ex{
i_{IJ}: \Ob_{IJ}\ar@{^{(}->}[r] &
\KGL(E,F) &\text{and}&
i_{r,s}: \Ob_{r,s}\ar@{^{(}->}[r] &
\KGL(E,F)
}
$$
denote the inclusion morphisms.

Let $\text{Fl}(E)$ and $\text{Fl}(F)$ denote the full flag manifolds 
associated to the vector spaces $E$ and $F$ respectively and
let $\Fl:=\text{Fl}(E)\times\text{Fl}(F)$.
For 
$a=(a_1,\dots,a_n)\in \Z^n$,\ 
$b=(b_1,\dots,b_n)\in \Z^n$
we define the invertible $\Oo_{\Fl}$-module
$$
\Oo_{\Fl}(a,b) := \Tensor_{i=1}^n(\E_i/\E_{i-1})^{\tensor a_i} \tensor
                    \Tensor_{i=1}^n(\F_i/\F_{i-1})^{\tensor b_i} 
\quad,
$$
where
$
0=\E_0\subset\E_1\subset\dots\subset\E_n=E\tensor\Oo_{\Fl}
$
and 
$
0=\F_0\subset\F_1\subset\dots\subset\F_n=F\tensor\Oo_{\Fl} 
$
are the two universal flags on $\Fl$.
The variety $\Fl$ is endowed with a canonical $\bG$-action and
the line bundles $\Oo_{\Fl}(a,b)$ come with a canonical $\bG$-linearization.

\begin{lemma}
\label{O_rs=Fl}
\label{pullback to O_rs}
For each pair $r,s\in[0,n]$ with $r+s=n$ we have a canonical isomorphism
$
\Ob_{r,s}\isomto\Fl
$,
which is compatible with the $\bG$-action on the two varieties.
Furthermore, we have a canonical 
isomorphism of $\bG$-linearized line bundles on $\Ob_{r,s}$:
$$
i_{r,s}^*
\left(
  \Tensor_{i=0}^{n-1}(M_i^{m_i}\tensor L_i^{l_i})\tensor 
  (\det E)^e\tensor(\det F)^d
\right)
=\Oo_{\Ob_{r,s}}(a,b)
\quad,
$$
where $\Oo_{\Ob_{r,s}}(a,b)$ is the line bundle corresponding
to $\Oo_{\Fl}(a,b)$ via the isomorphism $\Ob_{r,s}\isomto\Fl$ and
where $(a,b)\in\Z^n\times\Z^n$ is defined by
$$
a_i-e = -b_{n-i+1}+d =
\left\{
\begin{array}{ll}
l_{n-i+1}-l_{n-i} & \text{if $i\in[1,s]$}\\
m_{i-1}-m_i             & \text{if $i\in[s+1,n]$}
\end{array}
\right.
$$
(It is understood that $m_n=l_n=0$).
\end{lemma}

\begin{proof}
The first part of the lemma is a special case of \cite{kgl} Theorem 9.3:
Let
$$
U_i:=\left\{
\begin{array}{ll}
\E_i & \text{, if $0\leq i\leq s$} \\
\E_{i-1} & \text{, if $s+1\leq i\leq n+1$}
\end{array}
\right.
\quad
,
\quad
V_i:=\left\{
\begin{array}{ll}
\F_i & \text{, if $0\leq i\leq r$} \\
\F_{i-1} & \text{, if $r+1\leq i\leq n+1$}
\end{array}
\right.
$$
Then in the notation of loc. cit. we have
$$
\Ob_{r,s}=P_1\times_{\Fl}\dots\times_{\Fl}P_r\times_{\Fl}
          Q_s\times_{\Fl}\dots\times_{\Fl}Q_1\times_{\Fl}K'
$$
where 
$
P_p=\PGLb(V_{r-p+1}/V_{r-p},U_{s+p+1}/U_{s+p})
$,
$
Q_q=\PGLb(U_{s-q+1}/U_{s-q},V_{r+q+1}/V_{r+q})
$
and 
$
K'=\KGL(U_{s+1}/U_s,V_{r+1}/V_r)
$.
Since the bundles 
$V_{r-p+1}/V_{r-p}$,
$U_{s+p+1}/U_{s+p}$,
$U_{s-q+1}/U_{s-q}$,
$V_{r+q+1}/V_{r+q}$
are of rank one and the bundles
$U_{s+1}/U_s$,
$V_{r+1}/V_r$
are of rank zero, it follows
that
$
P_p=Q_q=K'=\Fl
$
and therefore $\Ob_{r,s}=\Fl$.

The second part of the lemma follows from the poof of \cite{kgl}, 9.3:
In the notation of that proof we have
$$
i_{r,s}^*M_i=\M_i
=\left\{
\begin{array}{lll}
\Oo_{\Fl} & \text{, if $s\in[0,s-1]$} \\
\M_0^{(1)} & \text{, if $i=s$} \\
\M_0^{(i-s+1)}\tensor(\M_0^{(i-s)})^{\dual} & \text{, if $s+1\leq i\leq n-1$}
\end{array}
\right.
$$
where the line bundle $\M_0^{(p)}$ is ingredient of the bf-morphism
$$
\left(
\xymatrix@C4ex{
\E_1^{(p)}
\ar[rr]^0_{(\M_0^{(p)},\mu_0^{(p)}=0)}
& & 
\E_0^{(p)}
\ar@/_1.2pc/|{\tensor}[ll]
}
\right)
$$
between $\E_0^{(p)}=U_{s+p+1}/U_{s+p}$ and $\E_1^{(0)}=V_{r-p+1}/V_{r-p}$.
But this means that we have a canonical isomorphism of line bundles
$$
\M_0^{(p)}=(U_{s+p+1}/U_{s+p})\tensor(V_{r-p+1}/V_{r-p})^{\dual}
$$
and consequently
$$
i_{r,s}^*M_i=
(U_{i+2}/U_{i+1})
\tensor
(U_{i+1}/U_{i})^{\dual}
\tensor
(V_{n-i+1}/V_{n-i})
\tensor
(V_{n-i}/V_{n-i-1})^{\dual}
$$
for $i=s,\dots,n-1$.
Analogously we have
$
i_{r,s}^*L_i=\Oo_{\Fl}
$
for $i=0,\dots,r-1$ and
$$
i_{r,s}^*L_i=
(V_{i+2}/V_{i+1})
\tensor
(V_{i+1}/V_{i})^{\dual}
\tensor
(U_{n-i+1}/U_{n-i})
\tensor
(U_{n-i}/U_{n-i-1})^{\dual}
$$
for $i=r,\dots,n-1$.
The stated formula follows from this together with the fact that
we have $\Oo_{\Fl}((1,\dots,1),(0,\dots,0))=\det(E)\tensor\Oo_{\Fl}$,
and
$\Oo_{\Fl}((0,\dots,0),(1,\dots,1))=\det(F)\tensor\Oo_{\Fl}$.
\end{proof} 

\begin{proposition}
\label{BBW}
Let $a=(a_1,\dots,a_n)$ and $b=(b_1,\dots,b_n)$ be two elements in $\Z^n$.
Then 
$
H^0(\Fl,\Oo_{\Fl}(a,b))\neq 0
$ 
if and only if $a$, $b$ are increasing,
i.e. if 
$a_1\leq a_2\leq \dots\leq a_n$ 
and
$b_1\leq b_2\leq \dots\leq b_n$.
The association
$$
(a,b)\mapsto H^0(\Fl,\Oo_{\Fl}(a,b))
$$
establishes a bijection between the set of 
all increasing $a, b\in\Z^n$ 
and the set of simple $\bG$-modules.
Furthermore, $H^p(\Fl,\Oo_{\Fl}(a,b))=0$ for all $p\geq 2$ and all increasing
$a,b\in\Z^n$.
\end{proposition}

\begin{proof}
This is a special case of the Borel-Bott-Weil
theorem (cf. e.g. \cite{Jantzen} II. 5.5).
\end{proof}

\section{Statement of the theorem}
\label{section theorem}

We keep the notations introduced in section 
\ref{preliminary}.

\begin{definition}
\label{A}
Let $L$ be a line bundle on $\KGL(E,F)$ of the form
$$
L=\Tensor_{i=0}^{n-1}(M_i^{m_i}\tensor L_i^{l_i})\tensor
  (\det E)^e\tensor(\det F)^d
\quad.
$$ 
Let $I,J\subseteq[0,n-1]$ and let $i_1:=\min(I)$, $j_1:=\min(J)$
where it is understood that $\min(\emptyset)=n$. Assume $i_1+j_1\geq n$.
We denote by $A_{IJ}(L)$ the set of all elements 
$(a,b)\in \Z^n\times\Z^n$, which have
the following properties:
\begin{enumerate}
\item
$a_1\leq a_2\leq\dots\leq a_n$ 
\item
$
\sum_{j=i+1}^{n}(a_j-e)\leq m_i
$
for all $i\in [n-j_1,n-1]$
and equality holds for $i\in I$.
\item
$
\sum_{j=1}^{n-i}(a_j-e)\geq -l_i
$
for all $i\in[n-i_1,n-1]$
and equality holds for $i\in J$.
\item
For all $i\in[1,n]$ the 
equality $a_i-e=-b_{n-i+1}+d$ holds.
\end{enumerate}
For abbreviation we denote by $A(L)$ the set $A_{\emptyset,\emptyset}(L)$.
\end{definition}

\begin{remark}
Notice that for $r,s\in[0,n]$ with $r+s=n$ the set 
$
A_{[s,n-1],[r,n-1]}(L)
$
contains at most the single element $(a,b)$ defined in \ref{pullback to O_rs}.
\end{remark}

\begin{theorem}
\label{theorem}
Let $L$ be a line bundle on $\KGL(E,F)$ of the form
$$
L=\Tensor_{i=0}^{n-1}(M_i^{m_i}\tensor L_i^{l_i})\tensor
  (\det E)^e\tensor(\det F)^d
$$
and let $I,J\subseteq[0,n-1]$ be subsets with $\min(I)+\min(J)\geq n$.
Then the following holds:

1.
The $\bG$-module $H^0(\Ob_{IJ},i_{IJ}^*L)$ 
comes with a canonical decomposition as follows:
$$
H^0(\Ob_{IJ},i_{IJ}^*L)= 
\Oplus_{(a,b)\in A_{IJ}(L)}H^0(\Fl,\Oo(a,b))
\quad.
$$

2. This decomposition is compatible with restriction in the
sense that the following is a commutative diagram of $\bG$-modules:
$$
\xymatrix{
H^0(\KGL,L) \ar@{=}[d] \ar[rr]^{\text{Res}} & &
H^0(\Ob_{I,J},i_{I,J}^*L) \ar@{=}[d] \\
\text{$
\underset{(a,b)\in A(L)}{\Oplus}\!\!\!\! 
H^0(\Fl,\Oo_{\Fl}(a,b))
$} 
\ar@<1ex>@{->>}[r] &
\text{$
\!\!\!\!\! 
\underset{(a,b)\in A(L)\cap A_{I,J}(L)}{\Oplus}\!\!\!\!\!\!\! 
H^0(\Fl,\Oo_{\Fl}(a,b))
$} 
\ar@<1ex>@{^(->}[r] &
 \text{$
\!\!\!\!
\underset{(a,b)\in A_{I,J}(L)}{\Oplus}\!\!\!\!\! 
H^0(\Fl,\Oo_{\Fl}(a,b))
$}
}
$$
where the lower arrows are the canonical projection and inclusion morphisms
induced by the inclusions $A(L)\cap A_{I,J}(L)\subseteq A(L)$ and
$A(L)\cap A_{I,J}(L)\subseteq A_{I,J}(L)$ respectively. 

3. Let 
$$
L'=\Tensor_{i=0}^{n-1}(M_i^{m'_i}\tensor L_i^{l'_i})\tensor
  f^*(\det E)^e\tensor f^*(\det F)^d
\quad,
$$
where $m'_i\leq m_i$ and $l'_j\leq l_j$ and equality holds,
if $i\in I$ and $j\in J$ respectively.
Then we have a commutative diagram of $\bG$-modules as follows:
$$
\xymatrix{
H^0(\Ob_{I,J},i_{I,J}^*L') 
\ar@{^(->}[r]^{\tensor\mu^{m-m'}\tensor\lambda^{l-l'}} 
\ar@{=}[d]
&
H^0(\Ob_{I,J},i_{I,J}^*L) 
\ar@{=}[d]
\\
\text{$
\underset{(a,b)\in A_{I,J}(L')}{\Oplus}H^0(\Fl,\Oo_{\Fl}(a,b))
$}
\ar@<1ex>@{^(->}[r]
&
\text{$
\underset{(a,b)\in A_{I,J}(L)}{\Oplus}H^0(\Fl,\Oo_{\Fl}(a,b))
$}
}
$$
where the upper horizontal arrow is induced by the section
$$
\left.
\left(
\mu_0^{m_0-m'_0}\tensor\dots\tensor
\mu_{n-1}^{m_{n-1}-m'_{n-1}}\tensor
\lambda_0^{l_0-l'_0}\tensor\dots\tensor
\lambda_{n-1}^{l_{n-1}-l'_{n-1}}
\right)\right|_{\Ob_{I,J}}
$$
of $i_{I,J}^*(L\tensor (L')^{-1})$
and the lower horizontal arrow is induced by the inclusion
$A_{I,J}(L')\subseteq A_{I,J}(L)$.
\end{theorem}

\section{Proof of the theorem}
\label{proof}

We fix a basis $(v_1,\dots,v_n)$ for $E$ and $(w_1,\dots,w_n)$ for $F$.
Let $B_1\subseteq\GL(E)$ and $B_2\subseteq\GL(F)$ 
be the Borel subgroups consisting of linear automorphisms fixing the flags 
$$
\text{
$
\{0\}\subset\langle v_1\rangle\subset
            \langle v_1,v_{2}\rangle\subset\dots\subset E
$
\qquad and \qquad
$
\{0\}\subset\langle w_n\rangle\subset
            \langle w_n,w_{n-1}\rangle\subset\dots\subset F
$
}
$$
respectively. 
Let $U_1\subset B_1$ and $U_2\subset B_2$ be the maximal
unipotent subgroups of $B_1$ and $B_2$ respectively.
Then $\bB:=B_1\times B_2$ is a Borel subgroup of $\bG$ and 
$\bU:=U_1\times U_2$
is its maximal unipotent subgroup.

Let $V:=\bU\times\Aa^n$ and 
let $\xi_1,\dots,\xi_n\in H^0(V,\Oo_V)$ be the pull back of the
coordinate functions on $\Aa^n$.
Let $V^o\subset V$ be the maximal open subset where all the $\xi_i$ are
invertible.
For every pair $r,s\in[0,n]$ with $r+s=n$ we have a morphism  
$
j_{r,s}^o: V^o\to\Isom(E,F)
$,
which on $R$-valued points ($R$ a $k$-algebra) is defined by
$$
j_{r,s}^o(x,y,z)=y\comp \zeta_r(z) \comp x
$$
for $x\in U_1(R)$, $y\in U_2(R)$, $z=(z_1,\dots,z_n)\in (R^{\times})^n$,
where  $\zeta_r(z):E\tensor_k R\isomto F\tensor_k R$ is the isomorphism 
defined with respect to the given basis of $E$ and $F$ by the diagonal
matrix $\diag(\zeta_{r,1}(z),\dots,\zeta_{r,n}(z))$ whose entries
are
$$
\zeta_{r,i}(z)=
\left\{
\begin{array}{ll}
z_{i}^{-1}\dots z_r^{-1} & \text{for $i\in[1,r]$} \\
z_{r+1}\dots z_i & \text{for $i\in[r+1,n]$}
\end{array}
\right.
$$

It follows from \ref{X(ell)} that the morphism $j^{(r,s)}$
extends to an open immersion
$$
j^{(r,s)}: V\To \KGL(E,F)
$$
whose image is the open affine subscheme $X(r)$.
We know from \ref{B-invariant} that $X(r)$ is $\bB$-invariant;
therefore the immersion 
$j^{(r,s)}$
induces a $\bB$-action on $V$.
Explicitly, on $R$-valued points
this action is given by
$$
b\cdot_{r}(x,y,z):=(\rho x \rho^{-1}u_1^{-1}, u_2\tau y\tau^{-1}, z')
\quad,
$$
where $b=(u_1\rho,u_2\tau)\in\bB(R)$,
$u_i\in U_i(R)$, 
$\rho=\diag(\rho_1,\dots,\rho_n)$ and 
$\tau=\diag(\tau_1,\dots,\tau_n)$ 
are $R$-valued points of 
the maximal torus of $B_1$ and $B_2$ respectively and
$z'=(z'_1,\dots,z'_n)\in R^n$  is defined by
$$
z'_i:=\left\{
\begin{array}{llll}
\rho_{i+1}^{-1}\rho_i\tau_{i+1}\tau_i^{-1}z_i & \text{if $i\in[1,r-1]$} \\
\rho_r\tau_r^{-1}z_r & \text{if $i=r$} \\
\rho_{r+1}^{-1}\tau_{r+1}z_{r+1} & \text{if $i=r+1$} \\
\rho_{i}^{-1}\rho_{i-1}\tau_i\tau_{i-1}^{-1}z_i & \text{if $i\in[r+2,n]$} 
\end{array}
\right.
$$

Let $I,J$ be two subsets of $[0,n-1]$ and assume $\min(I)+\min(J)\geq n$.
By this assumption there exist $r,s\in[0,n]$ with $r+s=n$ and 
$I\subseteq[s,n-1]$, $J\subseteq[r,n-1]$.
It is clear from \ref{X(ell)} that
the closed subscheme $V_{IJ}$ of $V$ defined by the
cartesian diagram
$$
\xymatrix{
V_{IJ} \ar@{^(->}[r] \ar[d]_{j^{(r,s)}_{IJ}} & V \ar[d]^{j^{(r,s)}} \\
\text{$\Ob_{IJ}$} \ar@{^(->}[r]_(.4){i_{IJ}} & \KGL(E,F)
}
$$
is cut out by the equations 
$\xi_{n-i}=0$ for $i\in I$ and $\xi_{i+1}=0$ for $i\in J$.

\begin{proposition}
\label{multiplicity 1}
Let $L$ be a $\bG$-linearized line bundle on $\KGL(E,F)$
and let $I,J\subseteq[0,n-1]$ with $\min(I)+\min(J)\geq n$.
Then for each simple $\bG$-module $W$ 
the $\bG$-module $H^0(\Ob_{IJ},i_{IJ}^*L)$ contains $W$ at most with 
multiplicity one as a submodule.
\end{proposition}

\begin{proof}
(Analogous to the proof of Lemma 8.2 in \cite{CP}).
Let $\s_1$ and $\s_2$ be two global sections of $i_{IJ}^*L$ which generate
$\bB$-invariant lines in $H^0(\Ob_{IJ},i_{IJ}^*L)$ on which $\bB$
operates by the same character. 
I claim that $\Ob_{IJ}$ contains a dense open $\bB$-orbit $\Omega$.
Indeed, let $z=(z_1,\dots,z_n)\in\Aa^n$, where $z_i=0$ if
$n-i\in I$ or $i-1\in J$, and $z_i=1$ else.
Then by the preceding formulae it follows easily that if we
choose $r,s\in[0,n]$ such that $r+s=n$, $I\subseteq [s,n-1]$, 
$J\subseteq[r,n-1]$, then
the image of the point $(1,z)\in V_{IJ}\subseteq\bU\times\Aa^n$ 
by the morphism $j^{(r,s)}_{IJ}$ is contained
in a dense open $\bB$-orbit in $\Ob_{IJ}$.
Therefore $\s_1/\s_2$ is a rational function on $\Ob_{IJ}$,
which is necessarily constant, since its restriction to $\Omega$ is constant.
\end{proof}

\begin{proposition}
\label{necessary}
Let $L=\Tensor_{i=0}^{n-1}(M_i^{m_i}\tensor L_i^{l_i})$ and let 
$I,J\subseteq[0,n-1]$
with $\min(I)+\min(J)\geq n$.
If the $\bG$-module $H^0(\Ob_{IJ},i_{IJ}^*L)$  contains an irreducible
$\bG$-module $W$ as a  submodule, then there exists an element
$(a,b)\in A_{IJ}(L)$ such that $W$ is isomorphic to the
$\bG$-module
$
H^0(\Fl,\Oo_{\Fl}(a,b))
$.
\end{proposition}

\begin{proof}
(Analogous to the proof of Proposition 8.2 in \cite{CP}).
Let $r,s\in[0,n]$ with $r+s=n$ and $I\subseteq[s,n-1]$, $J\subseteq[r,n-1]$.
Since $V_{IJ}$ is isomorphic to an affine space $\Aa^N$ for some $N$,
there exists a nowhere vanishing section $\s_0$ of the line bundle 
$L_{IJ}:=(j^{(r,s)}_{IJ})^*i_{IJ}^*L$ on $V_{IJ}$.  The group 
$\bB$ acts on $H^0(V_{IJ},L_{IJ})$ and for any $b\in \bB$ the 
section $b\cdot \s_0$ is again nowhere vanishing and thus a scalar
multiple of $\s_0$, since invertible functions on $V_{IJ}$ are constant. 
Therefore $\bB$ acts by a character on the line generated by $\s_0$.

Now let $\s\in W$ be a highest weight vector with respect to $\bB$.
Thus $\s$ is a global section of $i_{IJ}^*L$ which generates
a $\bB$-invariant line inside the space $H^0(\Ob_{IJ},i_{IJ}^*L)$. 
Let $\s_1$ be its pull back by $j^{(r,s)}_{IJ}:V_{IJ}\to\Ob_{IJ}$.
Then we have $\s_1=f\s_0$ for some regular function $f$ on $V_{IJ}$.
Clearly $f$ generates a $\bB$-invariant line in $H^0(V_{IJ},\Oo)$,
therefore $f$ is left unchanged by the action of the maximal unipotent 
subgroup $\bU$ and it follows
that $f$ must be a polynomial in the $\xi_i$, where $i\in[1,n]$, 
$n-i\notin I$, $i-1\notin J$. In fact $f$ must be a monomial in these
$\xi_i$, since otherwise the $\bB$-translates of $f$ would generate
a subspace of dimension $\geq 2$ of $H^0(V_{IJ},\Oo)$.

It follows from the above that there is a divisor 
$$
D=\sum_{i\in[0,n-1]} \beta_i Z_i \ \ + 
  \sum_{i\in[0,n-1]} \alpha_i Y_i
$$
on $\KGL(E,F)$, where $\beta_i\geq 0$, $\alpha_i\geq 0$ for all $i$ 
and $\beta_i=0$ if $i\in I$, $\alpha_i=0$ if $i\in J$ , such that
the pull back of $D$ to $\Ob_{IJ}\cap X(r)$ 
coincides with the restriction of the vanishing divisor of $\s$
to this open subscheme of $\Ob_{IJ}$.
Therefore there is a global section $\s'$ of $i_{IJ}^*L(-D)$ whose
image under the canonical map 
$i_{IJ}^*L(-D)\to i_{IJ}^*L$ is $\s$ and whose restriction to 
$\Ob_{IJ}\cap X(r)$ is nowhere vanishing. Since
the intersection $\Ob_{r,s}\cap X(r)$
is nonempty, it follows that the restriction of $\s'$ to
the closed subscheme $\Ob_{r,s}\subseteq\Ob_{IJ}$ is a nonzero 
section of $i_{r,s}^*L(-D)$.

Let 
$m'_i:=m_i-\beta_i$ and $l'_i:=l_i-\alpha_i$ and
let $(a,b)\in\Z^n\times\Z^n$ be defined by 
$$
a_i = -b_{n-i+1} =
\left\{
\begin{array}{ll}
l'_{n-i+1}-l'_{n-i} & \text{if $i\in[1,s]$}\\
m'_{i-1}-m'_i       & \text{if $i\in[s+1,n]$}
\end{array}
\right.
$$
with the convention that $m'_n=l'_n=0$.
By \ref{pullback to O_rs} we have $i_{r,s}^*L(-D)=\Oo_{\Ob_{r,s}}(a,b)$.

Consider the following diagram of $\bG$-modules:
$$
H^0(\Ob_{IJ},i_{IJ}^*L)\longleftarrow
H^0(\Ob_{IJ},i_{IJ}^*L(-D))\longrightarrow
H^0(\Ob_{r,s},\Oo_{\Ob_{r,s}}(a,b))
\quad.
$$
The left arrow is injective and maps $\s'$ to $\s$.
The right arrow maps $\s'$ to a non-zero element
and by \ref{O_rs=Fl} and \ref{BBW} the object on the right is a
simple $\bG$-module.
Therefore $H^0(\Ob_{r,s},\Oo_{\Ob_{r,s}}(a,b))$ is isomorphic
to $W$ as $\bG$-module.
Let us gather what we know about $(a,b)$:
\begin{enumerate}
\item
Since, as we have seen above, the line bundle $\Oo_{\Ob_{r,s}}(a,b)$
has a non-vanishing global section, it follows from \ref{BBW}
that $a_1\leq\dots\leq a_n$. 
\item
We have
$
\sum_{j=i+1}^{n}a_j=m'_i\leq m_i
$
for $i\in[s,n-1]$ and equality holds if $i\in I$.
\item
We have
$
\sum_{j=1}^{n-i}a_j=-l'_i\geq -l_i
$
for $i\in[r,n-1]$ and equality holds if $i\in J$.
\item
By definition, $a_i=-b_{n-i+1}$.
\end{enumerate}

Let $i_1:=\min(I)$ and $j_1:=\min(J)$.
In the above argument we can choose any $r,s$ with $r\in[n-i_1,j_1]$
and {\em a priori} $(a,b)$ depends on $r,s$
but by \ref{BBW} the fact that
$H^0(\Ob_{r,s},\Oo_{\Ob_{r,s}}(a,b))$ and $W$ are isomorphic
as $\bG$-modules determines $(a,b)$ which is therefore independent
of $r,s$.
It follows that the inequality in 2. holds for
all $i\in[n-j_1,n-1]$ and the inequality in 3. holds for all 
$i\in[n-i_1,n-1]$, i.e. we have $(a,b)\in A_{IJ}(L)$.
\end{proof}

Let $x_{ij}/x_{00}$ ($i,j\in[1,n]$) denote the coordinate functions
on $\GLn=\Isom(E,F)$ interpreted as rational functions on $\KGL(E,F)$.
For each integer $p\in[1,n]$ we define the rational function  
$d_p$ on $\KGLn$ as the determinant of the $p\times p$ sub-matrix
of $(x_{ij}/x_{00})_{i,j\in[1,n]}$ with indices in
$[1,p]\times[1,p]$, i.e. we set 
$ 
d_p:=\det_{[1,p][1,p]}(x_{ij}/x_{00})
$. For convenience we define $d_0:=1$.

\begin{lemma}
\label{divisor}
Let $a=(a_1,\dots,a_n)\in \Z^n$.
Then we have the following equality of divisors on $\KGL(E,F)$:
$$
\divi
\prod_{i=1}^{n}
\left(
\frac{d_i}{d_{i-1}}
\right)^{a_{n-i+1}}
=
\sum_{i=0}^{n-1}\left(-\sum_{j=i+1}^n a_j\right)Z_i +
\sum_{i=0}^{n-1}\left(\sum_{j=1}^{n-i} a_j\right)Y_i +
\sum_{i=1}^{n-1}(a_{n-i+1}-a_{n-i})\Delta_i
\quad,
$$
where $\Delta_i$ denotes the closure in $\KGL(E,F)$ of
the subscheme $\{d_i=0\}\subset\GL_n$.
Furthermore, for every $I,J\subseteq[0,n-1]$ with 
$\min(I)+\min(J)\geq n$ the closed subscheme
$\Ob_{IJ}$ is not contained in any of the $\Delta_i$.
\end{lemma}

\begin{proof}
The subvariety $\Delta_i$ is the locus of vanishing of the global section
$\det_{[1,i][1,i]}\Phi_{\univ}$ of the line bundle
$
\Tensor_{\nu=1}^i(\Tensor_{j=0}^{n-\nu}M_j \tensor 
                   \Tensor_{j=0}^{\nu-1}L_j^{-1})
$
and by \ref{det Phi} the complement of the union of all
$\Delta_i$ is precisely the union of the open sets $X(\ell)$
where $\ell$ runs through $[0,n]$.
In the notation introduced before \ref{X(ell)} we have 
$d_i/d_{i-1}=t_i/t_0$.
Now using \ref{X(ell)} a simple calculation shows that
for each $\ell$
the divisor of $\prod_{i=1}^n(t_i/t_0)^{a_{n-i+1}}$
on $X(\ell)$ is a linear combination of the restrictions
of the $Z_i$ and $Y_i$ with coefficients as given in the formula.

For the second part of the lemma
we choose $\ell$ such that $I\subseteq[n-\ell,n-1]$ and 
$J\subseteq[\ell,n-1]$
(which is always possible). Then the intersection of $\Ob_{IJ}$ with 
$X(\ell)$ is clearly nonempty; therefore $\Ob_{IJ}$ is not
contained in the complement of $X(\ell)$. In particular it is
not contained in any of the $\Delta_i$.
\end{proof}

We now come to the proof of Theorem \ref{theorem}.

\vspace{4mm}
{\em Proof of the first statement:}
Let $L$ and $I,J$ be as in the theorem and assume first that $e=d=0$.
If $A_{IJ}(L)$ is empty, then 
$H^0(\Ob_{IJ},i_{IJ}^*L)=(0)$
by \ref{necessary} and therefore the statement 1 of the theorem trivially 
holds in this case.

Assume $A_{IJ}(L)$ is non-empty,
let $(a,b)\in A_{IJ}(L)$ and let $r,s\in[0,n]$ with $r+s=n$ and
$I\subseteq[s,n-1]$, $J\subseteq[r,n-1]$.
Let $L':=\Tensor_{i=0}^{n-1}(M_i^{m'_i}\tensor L_i^{l'_i})$, where 
$m'_i:=\sum_{j=i+1}^n a_j$ and 
$l'_i:=-\sum_{j=1}^{n-i}a_j$.
From lemma \ref{divisor} it follows that there exists a global section
of $L'$ whose restriction to $\Ob_{r,s}$ is nonzero.
By \ref{pullback to O_rs} we have $i_{r,s}^*L'=\Oo_{\Ob_{r,s}}(a,b)$.
Therefore we have a non-vanishing restriction morphism 
$$
H^0(\Ob_{IJ},i_{IJ}^*L')\to H^0(\Ob_{r,s},\Oo_{\Ob_{r,s}}(a,b))
\quad.
$$
Together with \ref{BBW} and \ref{multiplicity 1} it follows 
that the $\bG$-module 
$H^0(\Ob_{IJ},i_{IJ}^*L')$ contains an irreducible 
submodule 
$W\isomorph H^0(\Ob_{r,s},\Oo_{\Ob_{r,s}}(a,b))$
exactly with multiplicity one.
In particular, by the above restriction morphism the submodule
$W \subseteq H^0(\Ob_{IJ},i_{IJ}^*L')$ 
is canonically identified with 
$H^0(\Ob_{r,s},\Oo_{\Ob_{r,s}}(a,b))$
which in turn is canonically isomorphic to 
$H^0(\Fl,\Oo_{\Fl}(a,b))$.

I claim that the identification of $W$ with 
$H^0(\Fl,\Oo_{\Fl}(a,b))$
is independent of the choice of the numbers $r,s$.
For this it is clearly sufficient to show
that the composite morphism
$$
H^0(\KGL,L')\to 
H^0(\Ob_{r,s},\Oo_{\Ob_{r,s}}(a,b))\isomto
H^0(\Fl,\Oo_{\Fl}(a,b))
\eqno(*)
$$
does not depend on $r,s$. 
This will be shown below. For the moment we assume this fact.

It is clear from the definition of $A_{IJ}(L)$ that
$$
\left.
\left(
 \mu_0^{m_0-m'_0}\tensor\dots\tensor
 \mu_{n-1}^{m_{n-1}-m'_{n-1}}\tensor
 \lambda_0^{l_0-l'_0}\tensor\dots
 \tensor\lambda_{n-1}^{l_{n-1}-l'_{n-1}}
\right)\right|_{\Ob_{IJ}}
$$
is a non-vanishing global section of $i_{IJ}^*(L\tensor (L')^{-1})$,
and therefore defines a canonical injective morphism
$$
H^0(\Ob_{IJ},i_{IJ}^*L')\to H^0(\Ob_{IJ},i_{IJ}^*L)
\quad.
$$
It follows that 
$H^0(\Ob_{IJ},i_{IJ}^*L)$
contains an irreducible $\bG$-submodule (the image of $W$),
which is canonically isomorphic to $H^0(\Fl,\Oo_{\Fl}(a,b))$.
This together with \ref{multiplicity 1} and \ref{necessary}
clearly implies statement 1 of the theorem in the case $e=d=0$.

For arbitrary $e,d$ the result is easily deduced from that special
case, the key observation being that we have a canonical isomorphism
$$
\Oo_{\Fl}(a,b)\tensor(\det E)^e\tensor(\det F)^d=
\Oo_{\Fl}(a+(e,\dots,e),b+(d,\dots,d))
$$
of $\bG$-linearized line bundles on $\Fl$.

\vspace{4mm}
{\em Independence of (r,s):}
It remains to be shown that the morphism $(*)$   
does not depend on $r,s$. In fact, since 
$H^0(\Fl,\Oo_{\Fl}(a,b))$
is a simple $\bG$-module, it suffices to produce a point $z$ in $\Fl$ 
and a global section
$\s$ of $L'$, whose image in $H^0(\Fl,\Oo_{\Fl}(a,b))$ evaluates to a 
nonzero element in in the fiber
$\Oo_{\Fl}(a,b)[z]$ of $\Oo_{\Fl}(a,b)$ at $z$,
which is independent of $r,s$.

By \ref{divisor}, the rational function 
$
\prod_{i=1}^{n}(d_i/d_{i-1})^{a_{n-i+1}}
=
\prod_{i=1}^n(t_i/t_0)^{a_{n-i+1}}
$ 
gives rise to a global section of $\Oo_{\KGL}(D)$, where
$$
D:=\sum_{i=0}^{n-1}(m'_iZ_i+l'_iY_i)
$$
Let $\s\in H^0(\KGL,L')$ be the element, which corresponds to this section
via the canonical isomorphism $\Oo_{\KGL}(D)\isomto L'$
and let $z\in\Fl$ be the point given by the pair of flags
$$
\text{
$
\{0\}\subset\langle v_n\rangle\subset\langle v_n,v_{n-1}\rangle
\subset\dots\subset E
$
\qquad and \qquad
$
\{0\}\subset\langle w_1\rangle\subset\langle w_1,w_2\rangle 
\subset\dots\subset F
$.
}
$$
Then the image of $\s$ by the morphism
\begin{eqnarray*}
&
H^0(\KGL,L') \to 
H^0(\Ob_{r,s},i_{r,s}^*L')\isomto
H^0(\Fl,\Oo_{\Fl}(a,b))\to
\Oo_{\Fl}(a,b)[z]=
&
\\
&
=
\Tensor_{i=1}^n
\left(
\frac{\langle v_n,\dots,v_{n-i+1}\rangle}
     {\langle v_n,\dots,v_{n-i+2}\rangle}
\right)^{\tensor a_i}
\tensor
\Tensor_{i=1}^n
\left(
\frac{\langle w_1,\dots,w_{i}\rangle}
     {\langle w_1,\dots,w_{i-1}\rangle}
\right)^{\tensor b_i}
=
\langle
\overset{n}{\underset{i=1}{\tensor}}
(v_i\tensor w_i^{-1})^{a_{n-i+1}}
\rangle
&
\end{eqnarray*}
is precisely the generator 
$
\tensor_{i=1}^{n}
(v_i\tensor w_i^{-1})^{a_{n-i+1}}
$
of the fiber of $\Oo_{\Fl}(a,b)$ at $z$.

Thus the image of $\s$ in $\Oo_{\Fl}(a,b)[z]$ does not depend on $r,s$
as was to be shown.

\vspace{4mm}
{\em Proof of the second statement:}
Let $(a,b)\in A(L)$ and let $(a',b')\in A_{I,J}(L)$. 
Consider the composite morphism
$$
H^0(\Fl,\Oo_{\Fl}(a,b))\To 
H^0(\KGL,L)\To H^0(\Ob_{I,J},i_{i,j}^*L)\To
H^0(\Fl,\Oo_{\Fl}(a',b'))
\eqno(\dagger)
$$
If $(a,b)\neq(a',b')$, then this morphism is clearly $0$, since
then the domain and the target are non-isomorphic simple $\bG$-modules.
If $(a,b)=(a',b')$, let 
$
L':=\Tensor_{i=0}^{n-1}(M_i^{m'_i}\tensor L_i^{l'_i})
    \tensor f^*(\det E)^e\tensor f^*(\det F)^d
$, 
where 
$m'_i:=\sum_{j=i+1}^n a_j$ and 
$l'_i:=-\sum_{j=1}^{n-i}a_j$,
let $r,s\in[0,n]$ with $I\subseteq[s,n-1]$ and $J\in[r,n-1]$
and
consider the commutative diagram
$$
\xymatrix@R=3ex{
H^0(\KGL,L) \ar[rr] & &
H^0(\Ob_{I,J},i_{I,J}^*L) \\
H^0(\KGL,L') \ar[rr] \ar@{^(->}[u] \ar@{->>}[d] & &
H^0(\Ob_{I,J},i_{I,J}^*L') \ar[d] \ar@{^(->}[u] \ar@{->>}[d] \\
H^0(\Ob_{r,s},i_{r,s}^*L') \ar@{=}[r] &
H^0(\Fl,\Oo_{\Fl}(a,b)) &
H^0(\Ob_{r,s},i_{r,s}^*L') \ar@{=}[l] &
}
$$
where the horizontal arrows are the restriction morphisms and 
the vertical arrows are defined as in the proof of the first 
statement of the theorem.
Since all $\bG$-modules in this diagram contain the simple submodule
$H^0(\Fl,\Oo_{\Fl}(a,b))$ with multiplicity one, it follows that
the morphism $(\dagger)$ has to be the identity in this case.

\vspace{4mm}
{\em Proof of the third statement:}
Let $(a',b')\in A_{I,J}(L')$ and let $(a,b)\in A_{I,J}(L)$. 
Consider the composite morphism
$$
H^0(\Fl,\Oo_{\Fl}(a',b'))\To 
H^0(\Ob_{I,J},i_{I,J}^*L')\To H^0(\Ob_{I,J},i_{i,j}^*L)\To
H^0(\Fl,\Oo_{\Fl}(a,b))
\eqno(\dagger\dagger)
$$
If $(a',b')\neq(a,b)$, this morphism vanishes by the
same argument as above.
If $(a',b')=(a,b)$, then the assertion that $(\dagger\dagger)$
is the identity morphism follows similarly as above from
the commutative diagram
$$
\xymatrix{
H^0(\Ob_{I,J},i_{I,J}^*L') 
\ar@{^(->}[rr]^{\tensor\mu^{m-m'}\tensor\lambda^{l-l'}} 
& &
H^0(\Ob_{I,J},i_{I,J}^*L) \\
H^0(\Ob_{I,J},i_{I,J}^*L'') \ar@{=}[rr] 
\ar@{^(->}[u]_{\tensor\mu^{m'-m''}\tensor\lambda^{l'-l''}} 
\ar@{->>}[d] & &
H^0(\Ob_{I,J},i_{I,J}^*L'') \ar[d] 
\ar@{^(->}[u]^{\tensor\mu^{m-m''}\tensor\lambda^{l-l''}} 
\ar@{->>}[d] \\
H^0(\Ob_{r,s},i_{r,s}^*L'') \ar@{=}[r] &
H^0(\Fl,\Oo_{\Fl}(a,b)) &
H^0(\Ob_{r,s},i_{r,s}^*L'') \ar@{=}[l] &
}
$$
where
$r,s\in[0,n]$ with $I\subseteq[s,n-1]$ and $J\in[r,n-1]$,
where 
$m''_i:=\sum_{j=i+1}^n a_j$, 
$l''_i:=-\sum_{j=1}^{n-i}a_j$ and
where
$
L'':=\Tensor_{i=0}^{n-1}(M_i^{m''_i}\tensor L_i^{l''_i})
    \tensor f^*(\det E)^e\tensor f^*(\det F)^d
$.

\section{Non-ampleness}
\label{remarks}

Since large parts of our proof of Theorem \ref{theorem}
are analogous to the proof of Theorem 8.3 in \cite{CP},
it is natural to ask whether also the analogue of the
ampleness result stated in Proposition 8.4 in \cite{CP} holds.
We will see below that the answer is negative.

Adopting the notation of \cite{CP} let $\Xb$ be the complete symmetric
variety associated to the data $(G,\sigma)$, where $G$ is a semi-simple
simply connected algebraic group over the complex numbers and $\sigma$
is a nontrivial involution on $G$. 
Let $S_1,\dots,S_\ell$ be the closures of the 1-codimensional
orbits of the natural action of $G$ on $\Xb$.
By \cite{CP}, 8.1 
the unique closed orbit $Y=\cap_{i=1}^\ell S_i$  
is isomorphic to $G/P$ for some parabolic
$P\subset G$ and the restriction of line bundles defines an injective
homomorphism $i^*:\Pic(\Xb)\to\Pic(Y)$. Proposition 8.4 in \cite{CP}
can be formulated as follows:

\vspace{2mm}
{\em Let $L$ be a line bundle on $\Pic(\Xb)$ such that $H^0(Y,i^*L)\neq 0$
and let $\omega_{\Xb}$ denote the dualizing line bundle on $\Xb$.
Then the line bundle $\omega_{\Xb}^{-1}(-S_1-\dots-S_\ell)\tensor L$ 
is ample.}

\vspace{2mm}
The following result shows that the analogue
of this Proposition in our context is false. 

\begin{proposition}
\label{not ample}
Let $\omega_{\KGL}:=\det\Omega^1_{\KGL(E,F)}$ 
denote the dualising line bundle on $\KGL(E,F)$.
Let $a_1\leq\dots\leq a_n$, $a_i\in\Z$ and let
$
L=\Tensor_{i=0}^{n-1}(M_i^{m_i}\tensor L_i^{l_i})
$,
where
$m_i:=\sum_{j=i+1}^n a_j$ and $l_i:=-\sum_{j=1}^{n-i} a_j$.
Then $H^0(\Ob_{r,s},i_{r,s}^*L)$ is nonzero for every
$r,s\in[0,n]$ with $r+s=n$, but neither the line bundle
$
\omega_{\KGL}^{-1}(-\sum_{i=0}^{n-1}(Z_i+Y_i))\tensor L
$
nor the line bundle $L$ itself is ample.
\end{proposition}

\begin{proof}
Using the inductive blowing up procedure which defines $\KGL(E,F)$
it is easy to see that 
$\omega_{\KGL}\isomorph\Tensor_{i=1}^{n-1}(M_i\tensor L_i)^{i(i-n)-1}$.

Let $e_1,\dots,e_n$ be the canonical basis of $N:=\Z^n$.
Let $\Delta$ be the smooth complete fan in $N_{\Q}=\Q^n$ whose
one-dimensional cones are generated by the vectors 
$\pm\sum_{i\in I}e_i$, where $I$ runs through the nonempty
subsets of $[1,n]$ and whose $n$-dimensional cones are
the sets
$
\sigma(\alpha,\ell):=
\{x\in\Q^n\ |\ x_{\alpha(1)}\leq\dots\leq x_{\alpha(\ell)}\leq 0\leq
               x_{\alpha(\ell+1)}\leq\dots\leq x_{\alpha(n)}\}
$
where $\alpha$ runs through the set of permutations of $[1,n]$
and $\ell$ runs through the set $[0,n]$.
Let $\KT$ be the smooth complete torus embedding associated to $\Delta$.

Let $\Tt$ be the torus embedding defined in \cite{kgl} p 563,
and let $Z_{i,\Tt}$, $Y_{i,\Tt}$ be the divisors on $\Tt$
defined there. The variety $\Tt$ can be identified with 
an open subscheme of $\KT$ and $\KT$ can be identified with
the closure in $\KGL(E,F)$ of a maximal torus in 
$\GL_n\isomorph\Isom(E,F)$ such that the restriction of the line
bundles $M_i$ and $L_i$ to $\Tt$ are $\Oo_{\Tt}(Z_{i,\Tt})$ and 
$\Oo_{\Tt}(Y_{i,\Tt})$ respectively. 

From \ref{pullback to O_rs} and \ref{BBW} it is immediate that
for any $r,s\in[0,n]$ with $r+s=n$ the restriction of $L$ to
$\Ob_{r,s}$ has non-vanishing global sections.
On the other hand, 
with the help of criterion \cite{Cox}, 3.1 it is easy to see 
that neither the restriction of $L$ nor that of   
$
\omega_{\KGL}^{-1}(-\sum_{i=0}^{n-1}(Z_i+Y_i))\tensor L
$
to the closed subscheme $\KT$ are ample.
Therefore the bundles $L$ and 
$
\omega_{\KGL}^{-1}(-\sum_{i=0}^{n-1}(Z_i+Y_i))\tensor L
$
cannot be ample.
\end{proof}

\end{document}